\documentclass[12pt]{amsart}
\usepackage{amscd,amssymb,amsthm}
\numberwithin{equation}{section}
\begin{document}

\title{Deformations of Special Lagrangian Submanifolds}
\author{Sema Salur}
 
\maketitle

\setcounter{section}{0}

\footnotesize

\vspace{.2in}

{\bf Abstract.} In [7], R.C.McLean showed that the moduli space of nearby submanifolds        of a smooth, compact, orientable special Lagrangian submanifold $L$ in a Calabi-Yau         manifold $X$ is a smooth manifold and its dimension is equal to the dimension of  $\mathcal{H}^1(L)$, the space of harmonic 1-forms on $L$. In this paper, we will show that the moduli space of all infinitesimal special Lagrangian deformations of $L$ in a symplectic manifold with non-integrable almost complex structure is also a smooth manifold of dimension $\mathcal{H}^1(L)$ . 

\vspace{.1in}

\normalsize
\section{Introduction}

The notion of special Lagrangian submanifold was first introduced by Harvey and Lawson [5] as an example of calibrated geometries. In the last few years, these have become more interesting after Strominger, Yau and Zaslow [9] proposed a geometric construction of mirror manifolds via special Lagrangian tori fibrations. According to their proposal, there should be a close connection between the deformation theory of special Lagrangian submanifolds and the mirror symmetry. Recently, special Lagrangian submanifolds have been studied by several other authors [1], [2], [3], [4], [6], [8].
\vspace{.1in}
 
In [7], McLean showed that the moduli space of nearby submanifolds        of a smooth, compact special Lagrangian submanifold $L$ in a Calabi-Yau         manifold $X$ is a smooth manifold and its dimension is equal to the dimension of  $\mathcal{H}^1(L)$, the space of harmonic 1-forms on $L$. In this paper, we will show that the moduli space of all infinitesimal special Lagrangian deformations of $L$ in a symplectic manifold with non-integrable almost complex structure is also a smooth manifold of dimension $\mathcal{H}^1(L)$. We will prove this by extending the parameter space of special Lagrangian
deformations, i.e, by using a modified definition of special Lagrangian submanifolds. Recall that a Lagrangian submanifold $L$ of a Calabi-Yau manifold is special Lagrangian if $Im(\xi)|_L\equiv 0$, where $\xi$ is a nowhere vanishing, closed, complex $(n,0)$-form. In our case, we will drop the assumption that $\xi$
is closed (i.e $d\xi\neq 0$) and introduce a new parameter $\theta$ for the deformations. Then the condition $Im(\xi)|_L\equiv 0$ will be replaced by $Im(e^{i\theta}\xi)|_L\equiv 0$ in the definition of special Lagrangian submanifolds.

\vspace{.1in}

In this paper, we will prove that the moduli space of all infinitesimal deformations of a smooth compact special Lagrangian submanifold $L$ in a symplectic manifold $X$ within the class of special Lagrangian submanifolds is a smooth manifold of dimension $\mathcal{H}^1(L)$.

\vspace{.1in} 

\section{Special Lagrangian Submanifolds}

In what follows, $X$ will denote a $2n$-dimensional symplectic manifold with symplectic 2-form $\omega$, an almost complex structure $J$ which is tamed by $\omega$, the compatible Riemannian metric $g$ and a nowhere vanishing complex valued $(n,0)$-form $\xi$$=$$\mu+i\beta$, where $\mu$ and $\beta$ are real valued $n$-forms. We say $\xi$ is {\em normalized} if the following condition holds:
 
\vspace{.1in}

\hspace{1in} $(-1)^{n(n-1)/2}(i/2)^n\xi\wedge\overline{\xi}$ = $\omega^n/n!$

\vspace{.1in}

So far, all the studies have focused on Calabi-Yau manifolds where this complex form is closed, but for our purposes we need a globally defined $(n,0)$-form which is not closed on $X$. There are many non-closed forms on a manifold, and one can easily construct them; for example by multiplying a given closed form with a non-constant function. The fact that $\xi$ is not closed implies that the associated almost complex structure $J$ on the tangent bundle $TX$ is non-integrable. In fact, the integrability of the almost complex structure determined by $\xi$ can be replaced by a weaker condition than $d\xi=0$ ([3]). 

\vspace{.1in}

In special Lagrangian calibrations, there is an additional term $e^{i\theta}$,
where for each fixed angle $\theta$ we have a corresponding form $e^{i\theta}\xi$ and its associated geometry. $\theta$ is called the {\em phase factor} of the calibration and this in fact will be our new parameter in the deformation of special Lagrangian submanifolds. In order to enlarge our parameter space we will allow $\theta$ to vary along the deformations. We will also assume that the initial value of $\theta$ is 0 to avoid the appearance of additional constants.

\vspace{.1in}

Taking the new parameter $\theta$ into consideration, one can slightly modify the definition of a special Lagrangian submanifold in a symplectic manifold.
                                                                   
\vspace{.1in}
                                   
{\bf Definition:} An $n$-dimensional submanifold $L\subseteq X$ is {\em special Lagrangian} if $L$ is Lagrangian (i.e. $\omega|_L\equiv 0$) and $Im(e^{i\theta}\xi)$ restricts to zero on $L$, for some $\theta\in {\bf R}$. Equivalently, $Re(e^{i\theta}\xi)$ restricts to be the volume form on $L$ with respect to the induced metric.

\vspace{.3in}

\section{Deformation Theory}

Under the given assumptions, we are ready to state our theorem:
                                                                               \vspace{.1in}

{\bf Theorem. }{\em   The moduli space of all infinitesimal deformations of a smooth, compact, orientable special Lagrangian submanifold $L$ in a symplectic manifold $X$ within the class of special Lagrangian submanifolds is a smooth manifold of dimension $\mathcal{H}^1(L)$ .}

\vspace{.1in}

{\bf Remark.} In [7], McLean proved the same theorem for Calabi-Yau manifolds, i.e, for  $d\xi=0$. This is the case where the almost complex structure is integrable.
 
 \vspace{.1in}

{\bf Proof of Theorem.}  For a small vector field $V$ and a scalar $\theta\in {\bf R}$, we define the deformation map as follows,

\vspace{.1in}     
$F: \Gamma(N(L))\times${\bf R}$\rightarrow \Omega^2(L)\bigoplus\Omega^n(L)$

\vspace{.1in}

$F(V,\theta)=((\exp_V)^*(-\omega), (\exp_V)^*(Im(e^{i\theta}\xi))$

\vspace{.1in}
 The deformation map $F$ is the restriction of $-\omega$ and $Im(e^{i\theta}\xi)$ to $L_V$ and then pulled back to $L$ via $(\exp_V)^*$ as in [7]. Here $N(L)$ denotes the normal bundle of $L$, $\Gamma(N(L))$ the space of sections of the normal bundle, and $\Omega^2(L)$, $\Omega^n(L)$ denote the differential $2$-forms and $n$-forms, respectively. Also, $\exp_V$ represents the exponential map which gives a diffeomorphism of $L$ onto its image $L_V$ in a neighbourhood of 0.

\vspace{.1in}
Recall that the normal bundle $N(L)$ of a special Lagrangian submanifold is isomorphic to the cotangent bundle $T^*(L)$. Thus, we have a natural identification of normal vector fields to $L$ with differential $1$-forms on $L$.
Furthermore, since $L$ is compact we can identify these normal vector fields with nearby submanifolds. Under these identifications, it is then easy to see that the kernel of $F$ will correspond to the special Lagrangian deformations.  

\vspace{.1in}
We compute the linearization of $F$ at (0,0),

\vspace{.1in}

$dF(0,0):\Gamma(N(L))\times${\bf R}$\rightarrow \Omega^2(L)\bigoplus\Omega^n(L)$ 
where

\vspace{.1in}

 $dF(0,0)(V,\theta)=\frac{\displaystyle\partial}{\displaystyle\partial{t}}F(tV,s\theta)|_{t=0, s=0} +\frac{\displaystyle\partial}{\displaystyle\partial{s}}F(tV,s\theta)|_{t=0, s=0} $

\vspace{.13in}
Therefore,
\vspace{.13in}

$\frac{\displaystyle\partial}{\displaystyle\partial{t}}F(tV,s\theta)|_{t=0, s=0}+\frac{\displaystyle\partial}{\displaystyle\partial{s}}F(tV,s\theta)|_{t=0, s=0}$
 
\vspace{.13in}

$=\frac{\displaystyle\partial}{\displaystyle\partial{t}}[\exp_{tV}^*(-\omega), \exp_{tV}^*(Im((\cos(s\theta)+i\sin(s\theta))(\mu+i\beta))]|_{t=0, s=0}$

\vspace{.13in}

\hspace{.1 in} $+\frac{\displaystyle\partial}{\displaystyle\partial{s}}[\exp_{tV}^*(-\omega), \exp_{tV}^*(Im((\cos(s\theta)+i\sin(s\theta))(\mu+i\beta))]|_{t=0, s=0}$

\vspace{.13in}

=$ [-({\it L}_{V}\omega)|_L, ({\it L}_{V}\mu)|_L\cdot \sin(s\theta)|_{s=0}
+  ({\it L}_{V}\beta)|_L\cdot \cos(s\theta)|_{s=0}$

\vspace{.13in}

\hspace{.1 in} + $((\exp_{tV}^*\mu)\cdot \cos(s\theta)\cdot \theta-(\exp_{tV}^*\beta)\cdot \sin(s\theta)\cdot \theta)|_{t=0, s=0}$]

\vspace{.2in}

 $=[-({\it L}_{V}\omega)|_L, {\it L}_{V}\beta|_L\cdot \cos(s\theta)|_{s=0} +((\exp_{tV}^*\mu)\cdot \cos(s\theta)\cdot \theta)|_{t=0, s=0}$]

\vspace{.1in}
Here ${\it L}_{V}$ represents the Lie derivative and one should notice that $\exp_{tV}^*\mu|_{t=0}$ is just the restriction of $\mu$ to $L$ which is equal to $1$ by our assumption that the initial value of $\theta$ is $0$.

\vspace{.1in}

Also, on a compact manifold $L$, top dimensional constant valued forms correspond to $\mathcal{H}^n(L)$, the space of harmonic $n$-forms on $L$ and there is a natural identification between the reals and harmonic $n$-forms. Therefore, $\theta$ will play the role of a harmonic $n$-form in our calculations.
  
\vspace{.1in}

Using the Cartan Formula, we get:

\vspace{.1in}

$=(-(i_Vd\omega +d(i_V\omega))|_L,(i_Vd\beta +d(i_V\beta))|_L + \theta)$
 
\vspace{.1in}

$=(-d(i_V\omega)|_L,(i_Vd\beta +d(i_V\beta))|_L +\theta)$

\vspace{.1in}

$=(dv, \zeta+d*v+\theta)$, where $\zeta=i_V(d\beta)|_L$

\vspace{.1in}

Here $i_V$ represents the interior derivative and $v$ is the dual $1$-form to the vector field $V$ with respect to the induced metric. For the details of local calculations of $d(i_V\omega)$ and $d(i_V\beta)$ see [7].

\vspace{.1in}

Hence $dF(0,0)(V,\theta) =(dv,\zeta+d*v+\theta)$.

\vspace{.1in}
Let $x_1, x_2,...,x_{n}$ and $x_1, x_2,...,x_{2n}$ be the local coordinates on $L$ and $X$, respectively. Then for any given normal vector field $V=(V_1\frac {\partial}{\partial x_{n+1}},...,V_n\frac {\partial}{\partial x_{2n}})$ to $L$ we can show that

\vspace{.1in}

$\zeta=i_V(d\beta)|_L=-n(V_1\cdot g_1+...+V_n\cdot g_n)$dvol where $g_i$ $(0<i\leq n)$ are combinations of coefficient functions in the connection-one forms. Since these are related to the induced metric one has the flexibility of making $g_i$ smaller by rescaling the metric on the ambient manifold. This fact will play an important role in proving the surjectivity of the linearized operator later on. 

\vspace{.1in}

One can also decompose the $n$-form $\zeta=da+d^*b+h_2$ by using Hodge Theory and because $\zeta$ is a top dimensional form on $L$, $\zeta$ will be closed and the
equation becomes  $dF(0,0)(V,\theta) =(dv,da+d*v+h_2+\theta)$ for some $(n-1)$-form $a$ and harmonic $n$-form $h_2$. 

\vspace{.1in}

The harmonic projection for $\zeta=-n(V_1.g_1+...+V_n.g_n)$dvol is $(\int_L{-n(V_1.g_1+...+V_n.g_n)}$dvol)dvol and therefore one can show that $da=-n(V_1.g_1+...+V_n.g_n)$dvol+$(n\int_L{(V_1.g_1+...+V_n.g_n)}$dvol)dvol and $h_2=(-n\int_L{(V_1.g_1+...+V_n.g_n)}$dvol)dvol.

\vspace{.1in}

{\bf Remark.} One should note that the differential forms $a$ and $h_2$ both depend on $V$ and therefore should be explored carefully in order to understand the deformations of special Lagrangian submanifolds.
  
\vspace{.1in}

After completing the space of differential forms with appropriate norms, we can consider $F$ as a smooth map from $C^{1,\alpha}(\Omega^1(L))\times{\bf R}$ to $C^{0,\alpha}(\Omega^2(L))$ and $C^{0,\alpha}(\Omega^n(L))$, where

\vspace{.1in}
$C^{k,\alpha}(\Omega)=\{ f\in C^k(\Omega) |$ $[D^{\gamma}f]_{\alpha,\Omega}<\infty , |\gamma|\leqq k \}$ and

\vspace{.16in}

$[f]_{\alpha,\Omega}= \mathop{Sup}\limits_{x,y\in \Omega,\; x\neq y } \frac{dist(f(x),f(y))}{(dist(x,y))^\alpha}$  in $\Omega$.
 
\vspace{.1in}

The Implicit Function Theorem says that $F^{-1}(0,0)$ is a manifold and its tangent space at $(0,0)$ can be identified with the kernel of $dF$. 

\vspace{.1in}

$(dv)\bigoplus(\zeta+d*v+\theta)=(0,0)$ implies 

$dv=0$ and $\zeta+d*v+\theta=da+d*v+h_2+\theta=0$. 

\vspace{.1in}

 The space of harmonic $n$-forms $\mathcal{H}^n(L)$, and the space of exact $n$-forms $d\Omega^{n-1}(L)$, on $L$ are orthogonal vector spaces by Hodge Theory. Therefore, $dv=0$ and $da+d*v+h_2+\theta=0$ is equivalent to $dv=0$ and $d*v+da=0$ and $h_2+\theta=0$.

\vspace{.1in}

One can see that the special Lagrangian deformations (the kernel of $dF$) can be identified with the $1$-forms on $L$ which satisfy the following equations:

\vspace{.1in}

(i) $dv=0$

(ii) $d*(v+\kappa(v))=0$ 

(iii) $h_2+\theta=0$.

\vspace{.1in}
Here, $\kappa(v)$ is a linear functional that depends on $v$ and $h_2$ is the harmonic part
of $\zeta$ which also depends on $v$. These equations can be formulated in a slightly different way in terms of decompositions of $v$ and $*a$.

\vspace{.1in}

If $v=dp+d^*q+h_1$ and $*a=dm+d^*n+h_3$ then we have 

\vspace{.1in}

(i) $dd^*q=0$

(ii) $\Delta(p\pm m)=0$

(iii) $h_2+\theta=0$.

\vspace{.1in}
This formulation of the solutions will help us to prove the surjectivity of the linearized operator without using $\kappa(v)$. 
\vspace{.1in}

{\bf Remark.} When $\theta=C$, the infinitesimal deformations of $\theta$ give no additional special Lagrangian deformations simply because there cannot be two different harmonic representatives in the same cohomology class. Therefore, one can obtain McLean's result by fixing $\theta=C$ along the deformations for some constant $C$ and since $d\beta|_L=0$ in the integrable case, $da=0$ and $h_2=0$. Hence the deformations correspond to 1-forms which satisfy the equations $dv=0$ and $d*v+\theta=0$.  

\vspace{.2in}

Next, we need to show that the deformation theory of special Lagrangian submanifolds is unobstructed. In order to use the implicit function theorem, we need to show that the linearized operator is surjective at $(0,0)$.
\vspace{.1in}

Recall that the deformation map,
\vspace{.1in}

$F:\Gamma(N(L))\times${\bf R}$\rightarrow\Omega^2(L)\bigoplus\Omega^n(L)$

\vspace{.1in}

is defined as follows:

\vspace{.1in}

$F(V,\theta)=((\exp_V)^*(-\omega), (\exp_V)^*(Im(e^{i\theta}\xi))$.

\vspace{.1in}
Even though $Im(e^{i\theta}\xi)$ is not closed on the ambient manifold $X$, the restriction of this differential form is a top dimensional form on $L$, and therefore it will be closed on $L$. On the other hand, $\omega$ is the symplectic form which is by definition closed on $X$. Therefore, the image of the deformation map $F$ lies in the closed $2$-forms and closed $n$-forms.

\vspace{.1in}

At this point we will investigate the surjectivity for $\omega$ and $Im(e^{i\theta}\xi)$ separately.

\vspace{.1in}

We have the following diagrams for $dF=dF_1\bigoplus dF_2$ with natural projection maps $proj_1$ and $proj_2$  :
\newcommand{\End}{\operatorname{End}}
\begin{equation*}
\begin{CD}
dF_1:\Gamma(N(L))@>d>>\Omega^2(L)@>proj_1>>d\Omega^1(L)
\end{CD}
\end{equation*}

and,

\begin{equation*}
\begin{CD}
dF_2:\Gamma(N(L))\times\bf R@>d*(1+\kappa)+\theta>>\Omega^n(L)@>proj_2>>d\Omega^{n-1}(L)\bigoplus\mathcal{H}^n(L)
\end{CD}
\end{equation*}

\vspace{.1in}

We will show that the maps $dF_1$ and $dF_2$ are onto $d\Omega^1(L)$ and $d\Omega^{n-1}(L)\bigoplus\mathcal{H}^n(L)$, respectively.
Therefore, for any given exact $2$-form $x$ and closed $n$-form $y=u+z$ in the image of the deformation map (here $u$ is the exact part and $z$ is the harmonic part of $y$), we need to show that there exists a 1-form $v$ and a constant $\theta$ that satisfy the equations,

\vspace{.1in}

(i) $dv=x$

(ii) $d*(v+\kappa(v))=u$

(iii)$h_2+\theta=z$.  
 
\vspace{.1in}

alternatively, we can solve the following equations for $p,q$ and $\theta$.

\vspace{.1in}

(i) $dd^*q=x$

(ii) $\Delta(p\pm m)=* u$  (Here, the star operator $*$ is defined on $L$)

(iii) $h_2+\theta=z$.

\vspace{.1in}

For (i), since $x$ is an exact $2$-form we can write $x=d(dr+d^*s+$harmonic form) by Hodge Theory. Then one can solve (i) for $q$ by setting $q=s$.

\vspace{.1in}

For (ii), since $\Delta m=d^*dm=d^**a=*d**a=\pm *da$,

$\Delta(p\pm m)=\Delta p \pm \Delta m= \Delta p \pm *da$ (here $a$ depends on $p$)

\vspace{.1in}

$= \Delta p \pm (-n(V_1.g_1+...+V_n.g_n)+(n\int_L{(V_1.g_1+...+V_n.g_n)}$dvol $))=*u$

\vspace{.1in}

Recall that for a given operator $\Delta + P:B\rightarrow B$ for some Banach space $B$, we can define its inverse by

\vspace{.1in}

$(\Delta +P)^{-1} = [\Delta ($Id$+\Delta ^{-1}P)]^{-1}$

\vspace{.1in}

\hspace{.6in}
$= [$Id$+\Delta ^{-1}P)]^{-1}\Delta ^{-1}$ 

\vspace{.1in}

\hspace{.6in}
$= [\sum_{k=0}^{\infty }{(-1)^k(\Delta ^{-1}P)^k]\Delta ^{-1}}$ 

\vspace{.1in}
and this is well defined if $||\Delta ^{-1}P||<1$.

\vspace{.1in}

The condition $||\Delta ^{-1}P||<1$ can be satisfied by taking $g_i$ very small. Since $g_i$ can be made smaller by rescaling the metric, one can make them small enough by multiplying the $(n,0)$-form $\xi$ with a small $\epsilon>0$. Recall that $(-1)^{n(n-1)/2}(i/2)^n\xi\wedge\overline{\xi}$ = $\omega^n/n!$. By making $\xi$ small we can make the symplectic $2$-form small and the associated metric $g$. Since $\epsilon$ is a real number the associated geometries of $\xi$ and $\epsilon.\xi$ will be equivalent. Since we can make $g_i$ small enough, we can invert the operator in (ii) and solve for $p$ for any given exact $n$-form $u$.

\vspace{.1in}

(iii) is straightforward.

\vspace{.1in}

The only thing remaining is to show that the image of the deformation map $F_1$ lies in $d\Omega^1(L)$ and the image of $F_2$ lies in $d\Omega^{n-1}(L)\bigoplus\mathcal{H}^n(L)$.

\vspace{.1in}

For $\omega$, we can follow the same argument as in [7]. Since $\exp_V:L$$\rightarrow $$X$ is homotopic to the inclusion $i:L$$\rightarrow$$X$, $\exp_V^*$ and $i^*$ induce the same map in cohomology. Thus, $[\exp_V^*(\omega)]=[i^*(\omega)]=[\omega|_L]=0$ . So the forms in the image of $F$ is cohomologous to zero. This is equivalent to saying that they are exact forms.
\vspace{.1in}

For  $Im(e^{i\theta}\xi$), we cannot follow the same process, because it is not a closed form on the ambient manifold $X$ and therefore does not represent a cohomology class. But by our construction of our deformation map, it is obvious that the image lies in $d\Omega^{n-1}(L)\bigoplus\mathcal{H}^n(L)$. 

\vspace{.1in}
One can find the dimension of this manifold by comparing the operators $d+*d^*(v)$ and $d+*d^*(v+\kappa(v))$. Since $\zeta=i_V(d\beta)|_L=-n(V_1.g_1+...+V_n.g_n)$\footnotesize dvol \normalsize it is easy to see that the extra term $*d^*(\kappa(v))$ contains no derivatives of $v$ and this implies that the linearized operators $d+*d^*(v)$ and $d+*d^*(v+\kappa(v))$ have the same leading term. Also it is known that the index of an elliptic operator is stable under lower order perturbations. Since the dimension of the kernel of $d+*d^*$ is $b_1(L)+1$ and the dimension of its cokernel is 1 as a map from $\Gamma(N(L))\times${\bf R}$\rightarrow d\Omega^1(L)\bigoplus d\Omega^{n-1}(L)\bigoplus \mathcal{H}^n(L)$, we can conclude that both the index of  $d+*d^*(v)$ and $d+*d^*(v+\kappa(v))$ are equal to $b_1(L)$. Hence the dimension of tangent space of special Lagrangian deformations in a symplectic manifold is also $b_1(L)$, the first Betti number of $L$. 

\vspace{.1in}

Therefore, $dF$ is surjective at $(0,0)$ and by infinite dimensional version of the implicit function theorem and elliptic regularity, the moduli space of all infinitesimal deformations of $L$ within the class of special Lagrangian submanifolds is a smooth manifold and has dimension $b_1(L)$. 

\vspace{.1in}

\small
{\em Acknowledgements.} This work was done when the author was visiting M.I.T. during the spring of 1999. Many thanks to the mathematics department at M.I.T. for their hospitality and support during the course of this work. Special thanks are to Gang Tian for all his help and encouragement. The author would also like to thank Jon Wolfson for introducing her to this subject. Also thanks to Mark Haskins for pointing out an error in the earlier version.

\normalsize

\begin{thebibliography}{[FP]}

\bibitem[1]{bf} Bryant, R.L. {\em Some examples of special Lagrangian Tori}, math.DG/9902076

\bibitem[2]{bf} Gross, M. {\em Special Lagrangian Fibrations I: Topology}, Integrable Systems and Algebraic Geometry (Kobe/Kyoto 1997), 156-193 World Scientific

\bibitem[3]{bf} Gross, M. {\em Special Lagrangian Fibrations II: Geometry}, alg-geom/9809072

\bibitem[4]{ab} Hitchin, N. {\em The moduli space of special Lagrangian submanifolds}, dg-ga/9711002

\bibitem[5]{b} Harvey, F.R. and Lawson, H.B. {\em Calibrated Geometries}, Acta. Math. {\bf 148} (1982), 47-157

\bibitem[6]{ab} Lu, P. {\em Special Lagrangian Tori on a Borcea-Voisin Threefold}, dg-ga/9902063

\bibitem[7]{bf} McLean, R.C. {\em Deformations of calibrated submanifolds}, Comm. Anal. Geom. {\bf 6} (1998), 705-747 

\bibitem[8]{bf} Merkulov, S.A. {\em The extended moduli space of special Lagrangian submanifolds}, math.AG/9806083

\bibitem[9]{bf} Strominger, A., Yau, S.T. and Zaslow, E., {\em Mirror Symmetry is T-Duality}, Nucl. Phys. {\bf B479} (1996), 243-259



\end{thebibliography}

\vspace{.2in} 

\small
Department of Mathematics,

Michigan State University, 

East Lansing, MI, 48824, USA

{\em E-mail address:} salur@math.msu.edu
\newpage
\end{document}